\documentclass[11pt]{article}
\usepackage{amsmath}
\usepackage{amssymb}
\usepackage{geometry}
\title{An Explicit Five-Variable Counterexample to the\\Generalized Vanishing Conjecture}
\author{Alexander Dvorsky\thanks{Department of Mathematics, University of Miami, Coral Gables, FL 33124, USA, dvorsky@math.miami.edu}}
\date{8/6/2026}

\begin{document}
\maketitle

\begin{abstract}
We give an explicit counterexample in five variables to the Generalized Vanishing Conjecture. The construction is motivated by the recent counterexample \cite{Long} to the Mathieu conjecture for $\mathrm{SU}(2)$. In the polynomial ring
\[
\mathbb{C}[a,b,c,d,t],
\]
set
\[
P=(t+c)(ad+bt), \qquad Q=c,
\]
and let
\[
\Lambda=\partial_t(\partial_a\partial_d-\partial_b\partial_c).
\]
We prove that
\[
\Lambda^m(P^m)=0 \qquad \text{for every } m\geq 1,
\]
whereas, for every $m\geq 2$,
\[
\Lambda^m(QP^m)=(-1)^m(m!)^2(m+1)!\,t\neq 0.
\]
Thus the Generalized Vanishing Conjecture fails in dimension 5. Both $P$ and the symbol of $\Lambda$ are homogeneous of degree 3.
\end{abstract}

\setcounter{section}{-1}
\section{Introduction}

Let
\[
R=\mathbb{C}[x_1,\ldots,x_n]
\]
and let $\Lambda\in\mathbb{C}[\partial_1,\ldots,\partial_n]$ be a constant-coefficient differential operator. The Generalized Vanishing Conjecture \cite[5.6]{EssenBook}, introduced in the context of Zhao's work on vanishing conjectures and Mathieu--Zhao spaces, asserts that
\[
\Lambda^m(f^m)=0 \qquad \text{for every } m\geq 1
\]
should imply
\[
\Lambda^m(gf^m)=0 \qquad \text{for all sufficiently large } m
\]
for every fixed polynomial $g\in R$.

We construct an explicit counterexample with $n=5$.

The starting point is Long's counterexample \cite{Long} to the Mathieu conjecture \cite[5.2]{EssenBook} for $\mathrm{SU}(2)$. Writing the coordinate functions on $\mathrm{SU}(2)$ as
\[
g=\begin{pmatrix}a&c\\ b&d\end{pmatrix},
\]
Long considers the regular functions
\[
F=(1+c)(ad+b), \qquad G=-c
\]
and proves that
\[
\int_{\mathrm{SU}(2)}F^m\,dg=0, \qquad
\int_{\mathrm{SU}(2)}GF^m\,dg=\frac{(-1)^{m-1}}{m+1}\neq 0
\]
for every $m\geq 1$.

The connection with differential operators comes from the Cayley operator
\[
\Delta=\partial_a\partial_d-\partial_b\partial_c.
\]
For a homogeneous polynomial $H\in\mathbb{C}[a,b,c,d]$ of degree $2r$, normalized Haar integration satisfies
\begin{equation}
\int_{\mathrm{SU}(2)}H(g)\,dg=\frac{1}{r!(r+1)!}\Delta^r H.
\end{equation}
Here $\Delta^rH$ is a constant. Formula (1) follows from the monomial identity
\begin{equation}
\int_{\mathrm{SU}(2)}a^\alpha b^\beta c^\gamma d^\delta\,dg
=(-1)^\beta\delta_{\alpha,\delta}\delta_{\beta,\gamma}
\frac{\alpha!\beta!}{(\alpha+\beta+1)!},
\end{equation}
which is the coordinate form of the $\mathrm{SU}(2)$ integration formula used by Mueger and Tuset \cite{MuegerTuset} and recalled in \cite{Long}.

The polynomial $F$ is not homogeneous. Introduce a homogenizing variable $t$ and set
\[
P=t^3F\!\left(\frac{a}{t},\frac{b}{t},\frac{c}{t},\frac{d}{t}\right)
=(t+c)(ad+bt).
\]
The operator
\[
\Lambda=\partial_t\Delta
\]
has order 3, equal to the degree of $P$. This explains the form of the construction. The verification below is entirely algebraic and does not require a general implication from the Mathieu conjecture to the Generalized Vanishing Conjecture.

During the preparation of this note, the author used OpenAI ChatGPT (model GPT-5.6 Sol, July 2026) as an interactive research assistant. ChatGPT assisted in testing the homogenized five-variable construction, generating computational verification code and locating relevant literature. All computations and references were then independently verified.

\section{The counterexample}

Set
\[
A=t+c, \qquad B=ad+bt,
\]
so that $P=AB$. Since all differential operators involved have constant coefficients,
\[
\Lambda^m=\partial_t^m\Delta^m.
\]
We first calculate $\Delta^m(P^m)$.

\subsection{Pure powers}

Expand
\begin{equation}
B^m=\sum_{j=0}^m\binom{m}{j}(ad)^{m-j}(bt)^j
\end{equation}
and
\begin{equation}
\Delta^m=\sum_{k=0}^m(-1)^k\binom{m}{k}
(\partial_a\partial_d)^{m-k}(\partial_b\partial_c)^k.
\end{equation}
When the $k$-th summand of (4) acts on $A^mB^m$, only the term with $j=k$ in (3) survives. Indeed, the $b$-derivatives require $j\geq k$, while the $a$- and $d$-derivatives require $j\leq k$.

After evaluating the derivatives, one obtains
\begin{align}
\Delta^m(P^m)
&=(m!)^2\sum_{k=0}^m(-1)^k\binom{m}{k}t^kA^{m-k}\notag\\
&=(m!)^2(A-t)^m\notag\\
&=(m!)^2c^m.
\tag{5}
\end{align}
The right-hand side is independent of $t$. Therefore
\begin{equation}
\Lambda^m(P^m)=\partial_t^m\Delta^m(P^m)=0 \qquad (m\geq 1).
\end{equation}
Thus $P$ satisfies the hypothesis of the Generalized Vanishing Conjecture.

\subsection{Mixed powers}

We now take $Q=c$. The same matching argument gives
\begin{equation}
\Delta^m(QP^m)=(m!)^2\sum_{k=0}^m\frac{(-1)^k}{k!}t^k\partial_c^k(cA^m).
\end{equation}
Since $A=t+c$, we may write
\[
cA^m=A^{m+1}-tA^m.
\]
Consequently,
\begin{equation}
\frac{1}{k!}\partial_c^k(cA^m)
=\binom{m+1}{k}A^{m+1-k}-t\binom{m}{k}A^{m-k}.
\end{equation}
Substituting (8) into (7) and applying the binomial theorem yields
\begin{align}
\Delta^m(QP^m)
&=(m!)^2\left[
\sum_{k=0}^m(-1)^k\binom{m+1}{k}t^kA^{m+1-k}
-t\sum_{k=0}^m(-1)^k\binom{m}{k}t^kA^{m-k}
\right]\notag\\
&=(m!)^2\left[c^{m+1}-tc^m+(-1)^mt^{m+1}\right].
\tag{9}
\end{align}
For $m\geq 2$, the first two terms in (9) are annihilated by $\partial_t^m$. Hence
\begin{align}
\Lambda^m(QP^m)
&=\partial_t^m\Delta^m(QP^m)\notag\\
&=(-1)^m(m!)^2(m+1)!\,t.
\tag{10}
\end{align}
This is nonzero for every $m\geq 2$. For completeness,
\[
\Lambda(QP)=-(c+2t)\neq 0
\]
when $m=1$.

We have therefore proved the following.

\medskip
\noindent\textbf{Theorem 1.1.} \textit{Let}
\[
P=(t+c)(ad+bt), \qquad Q=c,
\]
\textit{and}
\[
\Lambda=\partial_t(\partial_a\partial_d-\partial_b\partial_c)
\]
\textit{on $\mathbb{C}[a,b,c,d,t]$. Then}
\[
\Lambda^m(P^m)=0 \qquad \textit{for every } m\geq 1,
\]
\textit{but}
\[
\Lambda^m(QP^m)\neq 0 \qquad \textit{for every } m\geq 1.
\]
\textit{Consequently, the Generalized Vanishing Conjecture is false in dimension 5.}
\medskip

By adjoining unused variables, the conjecture also fails in every dimension $n\geq 5$.

\section{The Special Image Conjecture}

We record the corresponding explicit counterexample to the Special Image Conjecture \cite{Essen2010}, \cite[5.7]{EssenBook}. Let
\[
x=(a,b,c,d,t), \qquad \zeta=(\zeta_a,\zeta_b,\zeta_c,\zeta_d,\zeta_t),
\]
and put
\[
\mathcal{A}=\mathbb{C}[x,\zeta].
\]
Consider the subspace
\[
\mathcal{M}=\sum_{u\in\{a,b,c,d,t\}}(\partial_u-\zeta_u)\mathcal{A}.
\]
The Special Image Conjecture asserts that $\mathcal{M}$ is a Mathieu--Zhao subspace of $\mathcal{A}$. Thus, if $f^m\in\mathcal{M}$ for every $m\geq 1$, then for every fixed $g\in\mathcal{A}$ one should have $gf^m\in\mathcal{M}$ for all sufficiently large $m$.

Define the $\mathbb{C}$-linear map
\[
E:\mathcal{A}\longrightarrow\mathbb{C}[x]
\]
by
\[
E(\zeta^\alpha h(x))=\partial_x^\alpha h(x).
\]
It is standard that
\begin{equation}
\ker E=\mathcal{M}.
\end{equation}
See, for example, \cite{EssenBook}.

Let
\[
\lambda(\zeta)=\zeta_t(\zeta_a\zeta_d-\zeta_b\zeta_c),
\]
and retain the polynomials
\[
P=(t+c)(ad+bt), \qquad Q=c
\]
from the preceding section. Set
\[
f=\lambda(\zeta)P, \qquad g=Q.
\]
The constant-coefficient differential operator associated with $\lambda$ is
\[
\lambda(\partial)=\partial_t(\partial_a\partial_d-\partial_b\partial_c)=\Lambda.
\]
Consequently,
\begin{align}
E(f^m)&=E\bigl(\lambda(\zeta)^mP^m\bigr)\notag\\
&=\Lambda^m(P^m)\notag\\
&=0.
\tag{12}
\end{align}
for every $m\geq 1$. By (11),
\[
f^m\in\mathcal{M} \qquad (m\geq 1).
\]
On the other hand, for every $m\geq 2$,
\begin{align}
E(gf^m)&=E\bigl(Q\lambda(\zeta)^mP^m\bigr)\notag\\
&=\Lambda^m(QP^m)\notag\\
&=(-1)^m(m!)^2(m+1)!\,t\neq 0.
\tag{13}
\end{align}
Thus
\[
gf^m\notin\mathcal{M} \qquad (m\geq 2).
\]

\medskip
\noindent\textbf{Corollary 2.1.} \textit{The Special Image Conjecture is false in dimension 5. More explicitly, the subspace}
\[
\mathcal{M}=\sum_{u\in\{a,b,c,d,t\}}(\partial_u-\zeta_u)\mathbb{C}[x,\zeta]
\]
\textit{is not a Mathieu--Zhao subspace, as witnessed by}
\[
f=\zeta_t(\zeta_a\zeta_d-\zeta_b\zeta_c)(t+c)(ad+bt)
\]
\textit{and}
\[
g=c.
\]

By adjoining unused pairs of variables, the Special Image Conjecture also fails in every dimension $n\geq 5$.

\section{Remarks}

The conjectures relevant here form a one-way hierarchy. They were motivated initially by the Jacobian conjectures, and the implications are discussed in \cite[Chapter 5]{EssenBook}. Zhao's Special Image Conjecture implies the Generalized Vanishing Conjecture, and the Generalized Vanishing Conjecture in turn implies the Jacobian Conjecture. These implications are most naturally understood in their all-dimensional form: validity of the stronger conjecture in every finite dimension implies validity of the weaker one in every finite dimension. Their contrapositives therefore show that a counterexample to the Jacobian Conjecture forces the failure of both GVC and SIC somewhere, but they do not preserve dimension and need not produce an explicit counterexample. In particular, Alp\"oge's recent three-dimensional Keller map \cite{Sparkes} disproves the universal GVC and SIC, but the known implication chain does not convert it into a transparent counterexample to GVC(3), or even identify the smallest dimension in which the resulting vanishing statement fails.

Our example is homogeneous in a particularly strong sense:
\[
\deg P=3, \qquad \operatorname{ord}\Lambda=3,
\]
and the symbol of $\Lambda$ is
\[
\tau(\alpha\delta-\beta\gamma).
\]
The quadratic factor $\alpha\delta-\beta\gamma$ is irreducible, so $\Lambda$ is not a product of linear forms. This is consistent with the known positive results \cite{deBondt,EssenWillemsZhao} that the Generalized Vanishing Conjecture holds whenever the operator is a product of linear forms; in particular, it holds for every homogeneous constant-coefficient operator in two variables. As of now, GVC(3) and GVC(4) remain open.

\end{document}